\documentclass[conference]{IEEEtran}
\IEEEoverridecommandlockouts
\usepackage{cite}
\usepackage{amsmath,amssymb,amsfonts}
\usepackage{algorithmic}
\usepackage{graphicx}
\usepackage{textcomp}
\usepackage{xcolor}
\usepackage{url}
\usepackage{graphicx}
\usepackage{color}
\usepackage{placeins}
\usepackage{float}
\usepackage{tabularx,colortbl}
\usepackage{ifthen}

\newcommand{\vt}[1]{\boldsymbol{#1}} 
\newcommand{\vr}{\boldsymbol{r}} 
\newcommand{\mat}[1]{\mathbf{#1}}       
\renewcommand{\vec}[1]{\mathsfbfit{#1}}     
\newcommand{\op}[1]{\mathcal{#1}}           

\newcommand{\diff}{\,\mathrm{d}} 

\newcommand{\cald}{Calder\'on }

\newcommand{\im}{\mathrm{i}}  

\newcommand{\e}{\mathrm{e}}

\usepackage{graphicx}
\usepackage{subcaption}
\usepackage{amsfonts,amssymb,amsmath}
\usepackage{mathtools} 
\usepackage{hyperref} 
\usepackage{comment} 
\usepackage{soul}  
\usepackage{float} 
\usepackage{enumerate} 
\usepackage{siunitx} 
\usepackage{ulem} 
\usepackage{listings} 

\usepackage[utf8]{inputenc}
\usepackage{hyperref} 
\usepackage[T1]{fontenc}
\usepackage{subcaption}
\usepackage{float}
\usepackage{tikz}
\usetikzlibrary{arrows,matrix,shadows,positioning}
\usepackage{pgfplots}
\usepackage{pgfplotstable}
\usepackage{siunitx}
\usepackage{cleveref}
\pgfplotsset{compat=newest}
\usepackage{stmaryrd} 
\usepackage{amsthm}
\DeclareMathAlphabet{\mathbfsf}{\encodingdefault}{\sfdefault}{bx}{n}
\usepackage{bm}
\usepackage[OMLmathsfit]{isomath}

\usepackage{mathtools}

\usepackage{todonotes}
\input{title_singleColumn_URSI.tex}

\def\BibTeX{{\rm B\kern-.05em{\sc i\kern-.025em b}\kern-.08em
    T\kern-.1667em\lower.7ex\hbox{E}\kern-.125emX}}

\title{On some Operator Filtering Strategies Based on Suitably Modified Green's Functions}

\author[org1]{Matteo E. Masciocchi}
\author[org1]{Ermanno Citraro}
\author[org1]{Alexandre Dély}
\author[org1]{Lyes Rahmouni}
\author[org2]{Adrien Merlini}
\author[org1]{\newline Francesco P. Andriulli}

\address[org1]{Department of Electronics and Telecommunications, Politecnico di Torino, 10129 Turin, Italy}
\address[org2]{Microwave Department, IMT Atlantique, Brest, France}

\begin{document}

\newmaketitle
\begin{abstract}
Recent contributions showed the benefits of operator filtering for both preconditioning and fast solution strategies. 
While previous contributions leveraged laplacian-based filters, in this work we introduce and study a different approach leveraging the truncation of appropriately chosen spectral representations of  operators' kernels. In this contribution, the technique is applied to  the operators of the 2D TE- and TM-electric field integral equations (EFIE). We explore two different spectral representations for the 2D Green's function that lead to two distinct types of filtering of the EFIE operators.  Numerical results corroborate the effectiveness of the newly proposed approaches, also in the \cald preconditioned EFIE.
\end{abstract}

\section{Introduction}

The exploration of efficient solutions for complex electromagnetic problems is a pivotal research domain, underpinning advancements in numerous areas from communications to medical imaging. 
Integral equations are often used in this context, because they require the discretization of
scatterers' boundaries only. However, upon discretization via the boundary element method (BEM), they typically give rise to dense linear systems for which direct solutions are often unpractical.
One efficient approach to sidestep this difficulty is to couple fast solvers and iterative solvers to solve the system. However, both the solution's precision and the number of iterations required to reach it are strongly influenced by the spectral properties of the discretized integral operators.

Recently, operator filtering has been introduced \cite{merlini2022laplacian} as a way to
manipulate and correct the deleterious properties of electromagnetic operators'
spectra, while still being compatible with classical fast solution strategies. The first incarnations of operator filters rely on Laplacian manipulations to build spectral filters
that are subsequently multiplicatively applied to the integral operators of
interest. Operator filtering, relying on quasi-Helmholtz filters, has successfully been used to stabilize the Electric Field Integral Equations (EFIE) for 3D scattering 
in both the dense-discretization and low frequency regimes
\cite{merlini2022laplacian} and to enhance the compressibility of integral operators for building single-skeleton fast
direct solvers \cite{HenryFDS2022}.

In this work, we introduce a novel way to directly obtain filtered operators by truncating carefully chosen spectral representation of the operators' kernels. This means that the standard BEM discretization of the modified operators will directly yield matrices whose spectra correspond to a filtered version of the spectrum of the discretized original operators.
This article focuses on the operators involved in the TE and TM electric field integral equations (EFIE) applied to a 2D scatterer, although the proposed approach is extensible beyond this scenario, for which we present two types of filtering that rely on two different spectral representations of the 2D Green's function.
The properties and performance of the different filter approaches is analyzed in both the static and the dynamic case. 
The effectiveness of the new filtering scheme is further substantiated by suitably selected numerical results.

\section{Notation and Background}
\label{section:notation}
Consider a 2D scatterer modeled by a smooth curve $\gamma\in \mathbb{R}^2$, in a medium with wavenumber $k$ and impedance $\eta$, on which impinges an electromagnetic field $(\vt E^i, \vt H^i)$. 
We assume $\gamma$ lyes on the $xy$ plane, and denote with $t$ the physical quantities lying on the $xy$ plane, tangential to the scatterer.
Depending on the polarization of interest, two equations can be obtained to relate the incident electric field $\vt E^i = E_t \vt t + E_z \vt z$ to the induced surface current $\vt j = j_t \vt t + j_z \vt z$
\begin{align}
    \eta\im k\op{S}  j_z &=  E_\mathrm{z}\,,\\
    \eta\frac{1}{\im k}\op N  j_t &=  E_\mathrm{t}\,,
\end{align}
which are respectively the TM and TE EFIEs, and whose integral operators are
\begin{align}
    \label{eq:S}
    (\op S  j_z) (\vr) &\coloneqq  \int_{\gamma}  g(\vr, \vr') j_z(\vr') dr'\,,\\
    \label{eq:N}
    \left( \op N j_t \right)(\vr) &\coloneqq -\frac{\partial}{\partial n} \int_{\gamma} \frac{\partial}{\partial n'} g(\vr, \vr') j_t(\vr') \diff \vr'
\end{align}
with 
$g(\vr,\vr') \coloneqq -\im/4 H_0^{(2)}(k|\vr-\vr'|)$ if $k > 0$, or  
$g_0(\vr,\vr') \coloneqq -1/{2 \pi} \log(|\vr-\vr'|)$ if $k=0$ .
The currents $j_t$ and $j_z$ can be expanded with piecewise linear Lagrange interpolants $\{\varphi_i\}$ defined on a mesh of $\gamma$ composed of $N$ segments of uniform length $h$:
$j_t \approx\sum_{i=1}^N [\vec j_t]_i \varphi_i$ and $j_z\approx\sum_{i=1}^N [\vec j_z]_i \varphi_i$. 
Using Galerkin testing \cite{harrington_field_1993}, the discrete system can be obtained as
$\mat S \vec j_z = \vec E_{z}$
and
$\mat N \vec j_t = \vec E_{t}$,
where 
$\left[\mat S\right]_{ij} = \left<\varphi_i, \op S \varphi_j \right> $ , 
$\left[\mat N\right]_{ij} = \left<\varphi_i, \op N \varphi_j \right> $,
$\vec E_{t} = \left<\varphi_i, \frac{\im k}{\eta}E_t \right>$, 
and $\vec E_{z} = \left<\varphi_i, \frac{1}{\eta\im k}E_z \right>$.

\section{Operator Filtering via Green's Function Spectral Truncation}
\label{section:MGF}
In this section, we propose a new class of filtered operators, stemming from the truncation of spectral representations of the operators' kernel $g(\vr, \vr')$.
The scheme is as follows: (i) define the filtered Green's function $g^\alpha(\vr, \vr')$ where $\alpha$ indicates a filtering parameter (akin to a cutoff frequency) that will depend on the spectral representation chosen, (ii) define the filtered operators
\begin{align}
    \label{eq:Salpha}
    (\op S^\alpha  j_z) (\vr) &\coloneqq  \int_{\gamma}  g^\alpha(\vr, \vr') j_z(\vr') \diff r'\,, \\
    \label{eq:Nalpha}
    \left( \op N^\alpha j_t \right)(\vr) &\coloneqq -\frac{\partial}{\partial n} \int_{\gamma} \frac{\partial}{\partial n'} g^\alpha(\vr, \vr') j_t(\vr') \diff \vr' \,,
\end{align}
and (iii) use the boundary element method to obtain the matrices with the corresponding filtered spectra.

The first approach to filter $g(\vr - \vr')$ we will present, is to transform it into spectral domain in the sense of a multidimensional Fourier expansion, and back-transforming a truncated version obtaining  
 the following modified kernels in the static case
\begin{equation}
    \label{eq:2D_stat}
        g_0^\alpha(\vr, \vr') = -\frac{1}{2 \pi} \log(|\vr - \vr'|) - \frac{1}{2 \pi} \int_{s=\alpha}^{+\infty} \frac{J_0(s |\vr - \vr'|)}{s} \diff s\,,
\end{equation}
and in the dynamic case 
\begin{equation}
\label{eq:2D_dyn}
    g^\alpha(\vr, \vr') = -\frac{\im}{4} H_0^{(2)}(k|\vr - \vr'|) - \frac{1}{2 \pi} \int_{s=\alpha}^{+\infty} \frac{J_0(s |\vr - \vr'|) s}{s^2-k^2} \diff s
\end{equation}
respectively, where $\alpha > k$, $J_0$ is the $0$\textsuperscript{th} order Bessel function of the first kind and $H_0^{(2)}$ is the $0$\textsuperscript{th} order Hankel function of the second kind.
For the implementation of \eqref{eq:2D_dyn}, the computation is split in two parts: the singular part is handled by Taylor expansion; whereas the asymptotic regime is handled by a recursive extraction of terms by the expansion \cite[eq. 10.17.3]{olver2010nist}.

Another possible spectral expansion of $g(\vr, \vr')$ can be obtained leveraging Mehler–Sonine integrals \cite[eq. 10.9.12]{olver2010nist}, obtaining $Y_0 (x) = -\frac{2}{\pi} \int_{1}^{\infty} \frac{\cos (xt)}{\sqrt{(t^2-1)}} dt $, where $Y_0$ is the $0$\textsuperscript{th} order Bessel function of the second kind.
Using the identity $H_0^{(2)} (x) = J_0 (x) - \im  Y_0 (x) $, recalling the Green's Function definition, and truncating $ Y_0 (x)$, we obtain
\begin{equation}
\label{eq:1D}
    g^\alpha(\vr - \vr') = -\frac{\im}{4}J_0(k|\vr - \vr'|) - \frac{1}{2\pi} \int_{t=1}^{\alpha/k} \frac{\cos(k|\vr - \vr'| t)}{\sqrt{t^2 - 1}} \diff t\,.
\end{equation}
This form, however, is challenging to compute. An effective approach stems by rewriting the integral in \eqref{eq:1D} as
\begin{equation}
\label{eq:Ik1}
    I_k^{c1,c2} = \Re\left(\int_{c1}^{c2} \e^{\im kt}f(t)dt \right) .
\end{equation}
where $f:t\mapsto1/\sqrt{t^2-1}$. Now we introduce a change of variable which maps $\left(-1, 1\right)$ into $\left(c1, c2\right)$ carried out by a function 
$g:t\mapsto \left(t + 1\right) \frac{c2-c1}{2} + c1$.
Then, we expand $f \circ g$ as a linear combination of Legendre polynomials $P_n$, i.e., 
\begin{equation}
\label{eq:fPn}
    f(g(x)) = \sum_{n=0}^\infty a_n P_n(x), \ \,
    a_n = \frac{2n+1}{2} \int_{-1}^{1} f(g(x)) P_n(x)\diff x
\end{equation}
and so we can rewrite \eqref{eq:Ik1} as
\begin{equation}
    I_k^{c1,c2} = \Re\left(\frac{c2-c1}{2}\e^{\im \left(k' + k c1\right)} \sum_{n=0}^\infty a_n \int_{-1}^{1} \e^{\im k' t} P_n(t) dt \right)
\end{equation}
with $k' = (c2-c1)/2 k$ and the coefficients $a_n$  \eqref{eq:fPn}. 
Using the following identity from \cite{littlewood1976numerical}
\begin{equation}
    \int_{-1}^{1} P_n(t) \e^{\im k t} dt = (\im)^n \sqrt{\frac{2 \pi}{k}} J_{n+\frac{1}{2}}(k) 
\end{equation}
we finally obtain
\begin{equation}
    I_k^{c1,c2} = \Re\left(\frac{c2-c1}{2}\e^{\im \left(k' + kc1\right)} \sum_{n=0}^\infty a_n (\im)^n \sqrt{\frac{2 \pi}{k'}} J_{n+\frac{1}{2}}(k') \right) .
\end{equation}
In practice, the expansion above is used to compute the integral in $\left(2, \alpha/k\right)$, because in this way, due to the smoothness of the function over this interval, the Legendre expansion can be truncated with a low number of terms (and the Legendre coefficients $a_n$ can be precomputed for the integration interval of interest).
In the interval $\left(1,2\right)$ a different expansion is used: by integration by part of Eq. \eqref{eq:Ik1}, we obtain 
\begin{multline}
    I_k^{c1,c2} = \left[\cos(kt) \ln\left(t+\sqrt{t^2-1} \right) \right]_{t=c1}^{c2} + \\ k \int_{c1}^{c2} \sin(kt) \ln\left(t+\sqrt{t^2-1} \right) dt
\end{multline}
and the integral on the right hand side is computed using the Legendre expansion procedure. Because $\ln\left(t+\sqrt{t^2-1} \right)$ is not singular in $1$, the number of terms of the expansion is low, and the approach is efficient and accurate.

\section{Numerical results}
\label{section:numres}

\begin{figure}
    \centering
    \includegraphics[width=0.48\textwidth]{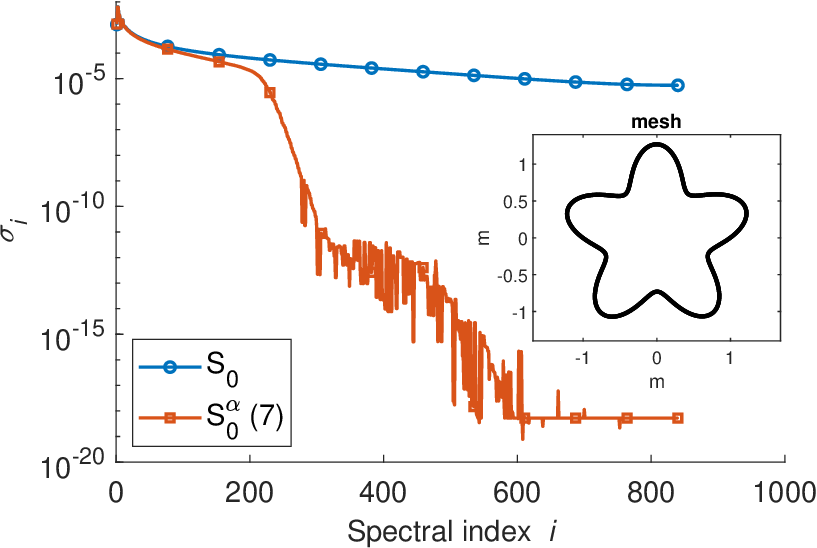}
    \caption{Singular values of $\mathbf{S}_0$ and singular values of $\mathbf{S}_0^\alpha$ using \eqref{eq:2D_stat}, ordered by the singular vectors of the Laplace-Beltrami operator, and reference mesh.}
    \label{fig:static}
\end{figure}

\begin{figure}
    \centering
    \includegraphics[width=0.48\textwidth]{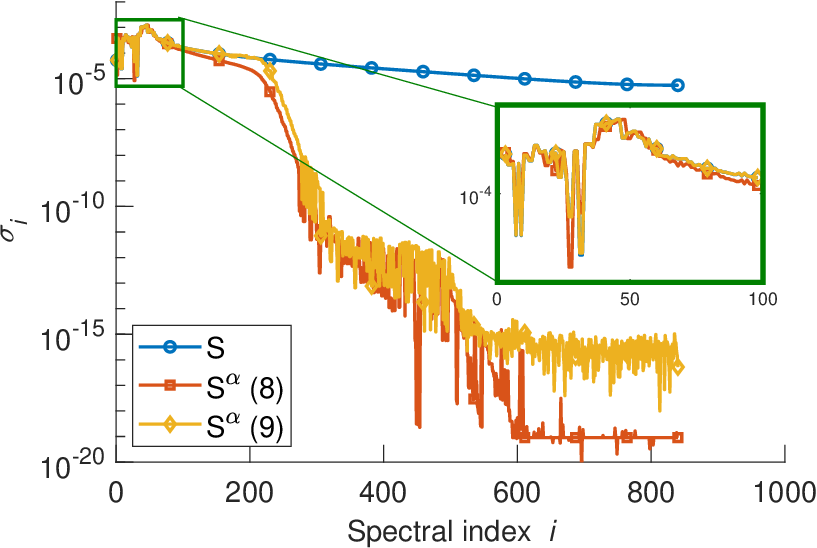}
    \caption{Singular values of $\mathbf{S}$ and singular values $\mathbf{S}^\alpha$ using \eqref{eq:2D_dyn} and \eqref{eq:1D}, ordered by the singular vectors of the Laplace-Beltrami operator.}
    \label{fig:dynS}
\end{figure}

\begin{figure}
    \centering
    \includegraphics[width=0.48\textwidth]{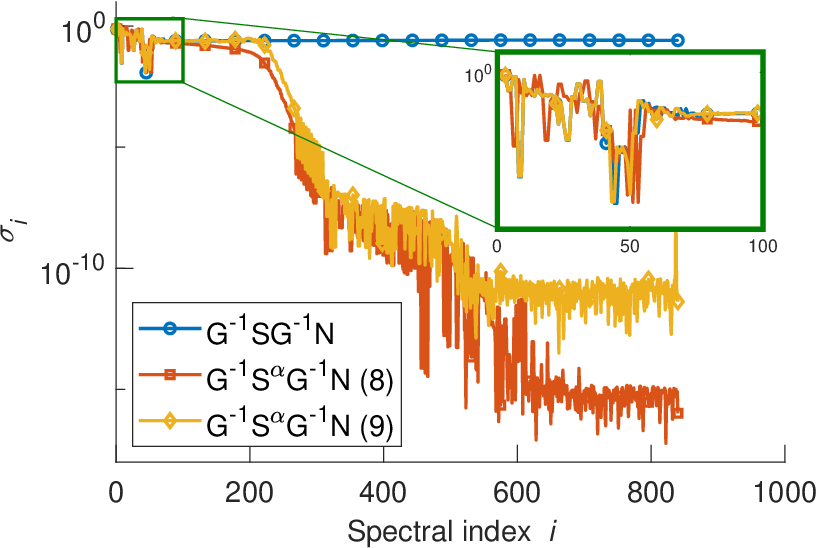}
    \caption{Singular values of $\mat G ^{-1} \mat S \mat G ^{-1} \mat N$ and singular values of $\mat G ^{-1} \mat S^\alpha \mat G ^{-1} \mat N$ where $\mathbf{S}^\alpha$ is computed using \eqref{eq:2D_dyn} and \eqref{eq:1D}, ordered by the singular vectors of the Laplace-Beltrami operator.}
    \label{fig:SN}
\end{figure}

To illustrate the effectiveness of the above formulations, numerical results for both the static and dynamic cases are provided below. 
All these operators were obtained on the scatterer shown in Fig.~\ref{fig:static}.
In the dynamic cases, the frequency is set to \SI{1}{\giga\hertz} with a mesh size of  $h = \lambda/30$.
Fig. \ref{fig:static} shows the effectiveness of \eqref{eq:2D_stat} on the discretized single layer operator $\mat S$ in static, whereas Fig. \ref{fig:dynS} show the performance of \eqref{eq:2D_dyn} and \eqref{eq:1D} on the same operator in the dynamic case. 
We finally show the effectiveness of the filtering procedure in the \cald precontitioned TE-EFIE \cite{adrianElectromagneticIntegralEquations2021} 
\begin{equation}
\label{eq:CaldTEdis}
    \mat G^{-1} \mat S 
    \mat G^{-1} \mat N \vec j_t =   
    \mat G^{-1} \mat S \mat G^{-1}\vec e_{t}\,,
\end{equation}
where $[\mat G]_{ij} \coloneqq \langle \varphi_i , \varphi_j \rangle$. In Figure~\ref{fig:SN},
in particular, we show the singular values of $ \mat G ^{-1} \mat S \mat G ^{-1} \mat N$ and singular values of $ \mat G ^{-1} \mat S^\alpha \mat G ^{-1} \mat N$ where $\mathbf{S}^\alpha$ is computed using \eqref{eq:2D_dyn} and \eqref{eq:1D}, ordered by the singular vectors of the Laplace-Beltrami operator. It is evident the effectiveness of the filtering procedure with the new kernels and thus their applicability in a \cald setting as further detailed in \cite{HenryFDS2022}.

\section*{Acknowledgment}

The work of this paper has received funding from the Horizon Europe Research and innovation programme under the EIC Pathfinder grant agreement n° 101046748 (project CEREBRO) and from the European Research Council (ERC) under the European Union’s Horizon 2020 research and innovation programme (grant agreement No 724846, project 321). 

\bibliographystyle{IEEEtran}
\bibliography{bibliography}

\end{document}